\documentclass[11pt,twoside]{article}
\usepackage[utf8]{inputenc}
\usepackage{kotex}
\usepackage{theorem}
\usepackage{epsfig,amsmath,amssymb}
\usepackage{color}
\allowdisplaybreaks[1]
\newtheorem{thm}{\bf Theorem}[section]

\newtheorem{lem}[thm]{\bf Lemma}

\newtheorem{rem}{\bf Remark}[section]
\numberwithin{equation}{section}
\newtheorem{exam}{\bf Example}[section]
\title{\bf Hyers-Ulam stability of isometries on bounded domains}

\author{Soon-Mo Jung}

\date{\small\it Mathematics Section,
                College of Science and Technology,
                Hongik University, 30016 Sejong,
                Republic of Korea \\
                \textit{E-mail:} {\tt smjung@hongik.ac.kr}}

\begin{document}

\markboth{Hyers-Ulam stability of isometries}{S.-M. Jung}
\maketitle

\begin{abstract}
More than 20 years after Fickett attempted to prove the
Hyers-Ulam stability of isometries defined on bounded subsets
of $\mathbb{R}^n$ in 1981, V\"{a}is\"{a}l\"{a} improved
Fickett's result significantly.
In this paper, we will improve Fickett's theorem by proving
the Hyers-Ulam stability of isometries defined on bounded
subsets of $\mathbb{R}^n$ using a more intuitive method
different from that used by V\"{a}is\"{a}l\"{a}.
\end{abstract}
\vspace{5mm}

\noindent
{\bf 2010 Mathematics Subject Classification}:
Primary 46C99; Secondary 39B82, 39B62, 46B04.
\vspace{5mm}

\noindent
{\bf Key Words}: Isometry; $\varepsilon$-isometry;
                 Hyers-Ulam stability;
                 bounded domain.

%%%%%%%%%%%%%%%%%%%%%%%%%%%
%%%   1. Introduction   %%%
%%%%%%%%%%%%%%%%%%%%%%%%%%%

\section{Introduction}

In 1940, Ulam gave a lecture at a mathematics club at the
University of Wisconsin introducing some important unsolved
problems.
Then, based on that lecture, he published a book 20 years
later (see \cite{U}).
A number of unresolved problems are introduced in this book,
among which the following question about Hyers-Ulam stability
of group homomorphism is closely related to the subject matter
of this paper:
\begin{quote}
\textit{Let $G_1$ be a group and let $G_2$ be a metric group
with the metric $d(\cdot, \cdot)$.
Given $\varepsilon > 0$, does there exist a $\delta > 0$ such
that if a function $h : G_1 \to G_2$ satisfies the inequality
$d(h(xy), h(x)h(y)) < \delta$ for all $x, y \in G_1$, then
there exists a homomorphism $H : G_1 \to G_2$ with
$d(h(x), H(x)) < \varepsilon$ for all $x \in G_1$?}
\end{quote}

In 1941, the following year, Hyers \cite{H} was able to
successfully solve Ulam's question about the approximately
additive functions, assuming that both $G_1$ and $G_2$ were
Banach spaces.
Indeed, he has proved that if a function $f : G_1 \to G_2$
satisfies the inequality
$\| f(x+y) - f(x) - f(y) \| \leq \varepsilon$ for some
$\varepsilon \geq 0$ and all $x, y \in G_1$, then there exists
an additive function $A : G_1 \to G_2$ such that
$\| f(x) - A(x) \| \leq \varepsilon$ for each $x \in G_1$.
In this case, the Cauchy additive equation,
$f(x+y) = f(x) + f(y)$, is said to have (or satisfy) the
Hyers-Ulam stability.
In the theorem of Hyers, the relevant additive function
$A : G_1 \to G_2$ is constructed from the given function $f$
by using the formula
$A(x) = \lim\limits_{n \to \infty} \frac{1}{2^n} f(2^n x)$.
This method is now called the \emph{direct method}.

We use the notations $(E, \| \cdot \|)$ and $(F, \| \cdot \|)$
to denote Hilbert spaces over $\mathbb{K}$, where $\mathbb{K}$
is either $\mathbb{R}$ or $\mathbb{C}$.
A mapping $f : E \to F$ is said to be an \emph{isometry} if
$f$ satisfies
\begin{eqnarray} \label{eq:2020-11-1-1}
\| f(x) - f(y) \| = \| x - y \|
\end{eqnarray}
for any $x, y \in E$.

By considering the definition of Hyers and Ulam \cite{hu1},
for each fixed $\varepsilon \geq 0$, a function $f : E \to F$
is said to be an \emph{$\varepsilon$-isometry} if $f$ satisfies
the inequality
\begin{eqnarray} \label{eq:2020-7-26-1}
\big| \| f(x) - f(y) \| - \| x - y \| \big| \leq \varepsilon
\end{eqnarray}
for any $x, y \in E$.
If there exists a positive constant $K$ depending only $E$ and
$F$ (independent of $f$ and $\varepsilon$) such that for each
$\varepsilon$-isometry $f : E \to F$, there is an isometry
$U : E \to F$ satisfying the inequality
$\| f(x) - U(x) \| \leq K \varepsilon$ for every $x \in E$,
then the functional equation (\ref{eq:2020-11-1-1}) is said to
have (or satisfy) the \emph{Hyers-Ulam stability}.

To the best of our knowledge, Hyers and Ulam were the first
mathematicians to study the Hyers-Ulam stability of isometries
(see \cite{hu1}).
Indeed, they were able to prove the stability of isometries
based on the properties of the inner product of Hilbert spaces:
For each surjective $\varepsilon$-isometry $f : E \to E$
satisfying $f(0) = 0$, there is a surjective isometry
$U : E \to E$ satisfying $\| f(x) - U(x) \| \leq 10 \varepsilon$
for every $x \in E$.
We encourage readers who want to read historically important
papers dealing with similar topics to look for the papers
\cite{dgbo1,dgbo2,rdbo,hu2}.

In 1978, Gruber \cite{gru} proved the following theorem:
Assume that $E$ and $F$ are real normed spaces, $f : E \to F$
is a surjective $\varepsilon$-isometry, and that $U : E \to F$
is an isometry satisfying $f(p) = U(p)$ for a $p \in E$.
If $\| f(x) - U(x) \| = o( \| x \| )$ as $\| x \| \to \infty$
uniformly, then $U$ is a surjective linear isometry and
$\| f(x) - U(x) \| \leq 5 \varepsilon$ for all $x \in E$.
In particular, if $f$ is continuous, then
$\| f(x) - U(x) \| \leq 3 \varepsilon$ for each $x \in E$.
This Gruber's result was further improved by Gevirtz \cite{gev}
and by Omladi\v{c} and \v{S}emrl \cite{os}.
There are many other papers related to the stability of
isometries, but it is regrettable that due to the restriction
of space, they cannot be quoted one by one.
Nevertheless, see
\cite{baker,bs,dol,fic,jung5,jung6,ls,ras1,ras,rs,sem,sko,
skof,swa,va} for more general information on the stability of
isometries and related topics.

The following Fickett's theorem is an important motive for this
paper (see \cite{fic}):

\begin{thm}[Fickett] \label{thm:fickett}
For a fixed integer $n \geq 2$, let $D$ be a bounded subset
of\/ $\mathbb{R}^n$ and let $\varepsilon > 0$ be given.
If a function $f : D \to \mathbb{R}^n$ satisfies the inequality
$(\ref{eq:2020-7-26-1})$ for all $x, y \in D$, then there
exists an isometry $U : D \to \mathbb{R}^n$ such that
\begin{eqnarray} \label{eq:fickett}
\| f(x) - U(x) \| \leq 27 \varepsilon^{1/2^n}
\end{eqnarray}
for each $x \in D$.
\end{thm}

The upper bound associated with inequality (\ref{eq:fickett})
in Fickett's theorem becomes very large for any sufficiently
small $\varepsilon$ in comparison to $\varepsilon$.
This is a big drawback of Fickett's theorem.
Thus, the work of further improving Fickett's theorem has to
be attractive.

After Fickett attempted to prove the Hyers-Ulam stability of
isometries defined on a bounded subset of $\mathbb{R}^n$ in
1981, several papers have been published steadily to improve
his result over the past 40 years.
However, most of the results were not very satisfactory
(see \cite{36,choi,jung5,jung6,jung2000}).
Fortunately, however, V\"{a}is\"{a}l\"{a} \cite{va} introduced
a new concept called the $c$-solar system and significantly
improved Fickett's result by proving the Hyers-Ulam stability
of isometries defined on bounded subsets of $\mathbb{R}^n$.

In this paper, we significantly improve Fickett's theorem by
using a more intuitive method that is different from
V\"{a}is\"{a}l\"{a}'s method.
Indeed, we prove the Hyers-Ulam stability of isometries defined
on bounded subsets of $\mathbb{R}^n$ for $n \geq 3$
(see Theorem \ref{thm:ndtheorem}).

%%%%%%%%%%%%%%%%%%%%%%%%%%%%%%%%%%%%%%%%%
%%%%%     2. Preliminaries     %%%%%%%%%%
%%%%%%%%%%%%%%%%%%%%%%%%%%%%%%%%%%%%%%%%%

%\section{Preliminaries}

%%%%%%%%%%%%%%%%%%%%%%%%%%%%%%%%%%%%%%%%%%%%%%%%%%%%%%%%%%%%
%%%%%     2. Real version of QR decomposition     %%%%%%%%%%
%%%%%%%%%%%%%%%%%%%%%%%%%%%%%%%%%%%%%%%%%%%%%%%%%%%%%%%%%%%%

\section{Real version of QR decomposition}

An \textit{orthogonal matrix} $\textbf{Q}$ is a real square
matrix whose columns and rows are orthonormal vectors.
In other words, a real square matrix $\textbf{Q}$ is orthogonal
if its transpose is equal to its inverse:
$\textbf{Q}^{\,tr} = \textbf{Q}^{-1}$, where $\textbf{Q}^{\,tr}$
and $\textbf{Q}^{-1}$ stand for the transpose and the inverse
of $\textbf{Q}$, respectively.
As a linear transformation, an orthogonal matrix preserves the
inner product of vectors, and therefore acts as an isometry of
Euclidean space.

Most papers and textbooks that mention \textit{QR decomposition}
only prove the complex version of QR decomposition
(see \cite[Theorem 6.3.7]{anton} or
\cite[Theorem 2.2 in \S 1]{stewart}).
However, not the complex version but the real version of QR
decomposition is required in this paper.
Nevertheless, since the proof of the real version is similar
to the proof of the complex version, we will omit the proof of
the real version here.
%Therefore, in the following theorem, we attempt to prove the
%real version of QR decomposition.
%By doing so, the completeness of this paper will be improved.

%%%%%%%%%%%%%%%%%%%%%%%%%%%%%%%%%%%%%%%%%%%%
%%%%%     QR decomposition     %%%%%%%%%%%%%
%%%%%%%%%%%%%%%%%%%%%%%%%%%%%%%%%%%%%%%%%%%%

\begin{thm}[QR decomposition] \label{thm:QR}
Every real square matrix {\rm $\textbf{A}$} can be decomposed
as {\rm $\textbf{A} = \textbf{Q} \textbf{R}$}, where
{\rm $\textbf{Q}$} is an orthogonal matrix and {\rm $\textbf{R}$}
is an upper triangular matrix whose elements are real numbers.
In particular, every diagonal element of {\rm $\textbf{R}$} is
nonnegative.
\end{thm}
\vspace*{3mm}

We can prove the following lemma using the real version of QR
decomposition (Theorem \ref{thm:QR}), and this lemma plays an
important role in achieving the final goal of this paper.
In practice, using this lemma, we can almost halve the number
of unknowns to consider in the main theorem.

%%%%%%%%%%%%%%%%%%%%%%%%%%%%%%%%%%
%%%%%     Lemma 2.3     %%%%%%%%%%
%%%%%%%%%%%%%%%%%%%%%%%%%%%%%%%%%%

\begin{lem} \label{lem:trans}
Let $\mathbb{R}^n$ be the $n$-dimensional Euclidean space for
a fixed integer $n > 0$.
Assume that $D$ is a subset of $\mathbb{R}^n$ with
$\{ e_1, e_2, \ldots, e_n \} \subset D$, where
$\{ e_1, e_2, \ldots, e_n \}$ is the standard basis for
$\mathbb{R}^n$.
If $f : D \to \mathbb{R}^n$ is a function and every $e_i$ is
written in column vector, then there exist an orthogonal matrix
{\rm $\textbf{P}$} and real numbers $e_{ij}^\prime$ for
$i, j \in \{ 1, 2, \ldots, n \}$ with $i \geq j$ such that
\begin{eqnarray*}
\emph{\textbf{P}} f(e_i)
= \left( e_{i1}^\prime, e_{i2}^\prime, \ldots, e_{ii}^\prime,
         0, \ldots, 0
  \right)^{tr}
\end{eqnarray*}
for every $i \in \{ 1, 2, \ldots, n \}$.
In particular, $e_{ii}^\prime \geq 0$ for all
$i \in \{ 1, 2, \ldots, n \}$.
\end{lem}
\vspace*{3mm}

\noindent
{\em Proof}.
Let $f(e_i) = ( e_{i1}, e_{i2}, \ldots, e_{in} )^{tr}$, written
in column vector, for any $i \in \{ 1, 2, \ldots, n \}$.
We now define a matrix $\textbf{A}$ by
\begin{eqnarray*}
\textbf{A}
= \big( f(e_1) ~ f(e_2) ~ \cdots ~ f(e_n) \big)
= \left( \begin{array}{cccc}
            e_{11} & e_{21} & \cdots & e_{n1} \\
            e_{12} & e_{22} & \cdots & e_{n2} \\
            \vdots & \vdots & \ddots & \vdots \\
            e_{1n} & e_{2n} & \ldots & e_{nn}
         \end{array}
  \right).
\end{eqnarray*}
By Theorem \ref{thm:QR}, there exist an orthogonal matrix
$\textbf{Q}$ and an upper triangular (real) matrix $\textbf{R}$
such that $\textbf{A} = \textbf{Q} \textbf{R}$.
Thus, we have
\begin{eqnarray} \label{eq:2021-1-27}
\textbf{R} = \textbf{Q}^{\,tr} \textbf{A}
= \left( \textbf{Q}^{\,tr} f(e_1) ~~
         \textbf{Q}^{\,tr} f(e_2) ~~ \cdots ~~
         \textbf{Q}^{\,tr} f(e_n)
  \right)
= \left(
  \begin{array}{cccc}
     e_{11}^\prime & e_{21}^\prime & \cdots & e_{n1}^\prime \\
     0             & e_{22}^\prime & \cdots & e_{n2}^\prime \\
     \vdots        & \vdots        & \ddots & \vdots \\
     0             & 0             & \ldots & e_{nn}^\prime
  \end{array}
  \right)
\end{eqnarray}
for some real numbers
$e_{11}^\prime, e_{21}^\prime, e_{22}^\prime, \ldots,
 e_{n1}^\prime, \ldots, e_{nn}^\prime$, where the diagonal
element $e_{ii}^\prime$ is nonnegative for all
$i \in \{ 1, 2, \ldots, n \}$.

Finally, the last two terms of (\ref{eq:2021-1-27}) are
compared to conclude as follows:
\begin{eqnarray*}
\left\{
\begin{array}{l}
   \textbf{Q}^{\,tr} f(e_1)
   = \left( e_{11}^\prime, 0, 0, \ldots, 0 \right)^{tr},
   \vspace*{2mm}\\
   \textbf{Q}^{\,tr} f(e_2)
   = \left( e_{21}^\prime, e_{22}^\prime, 0, \ldots, 0
     \right)^{tr},
   \vspace*{2mm}\\
   \hspace*{17.5mm} \vdots \vspace*{2mm}\\
   \textbf{Q}^{\,tr} f(e_n)
   = \left( e_{n1}^\prime, e_{n2}^\prime, e_{n3}^\prime,
            \ldots, e_{nn}^\prime
     \right)^{tr}.
\end{array}
\right.
\end{eqnarray*}
If we put $\textbf{P} = \textbf{Q}^{\,tr}$, then $\textbf{P}$
is also an orthogonal matrix.
\hspace*{\fill}$\Box$
\vspace*{3mm}

%%%%%%%%%%%%%%%%%%%%%%%%%%%%%%%%%%%%%%%%%%%%%%%
%%%%%     3. A preliminary theorem     %%%%%%%%
%%%%%%%%%%%%%%%%%%%%%%%%%%%%%%%%%%%%%%%%%%%%%%%

\section{A preliminary theorem}

Let $\{ e_1, e_2, \ldots, e_n \}$ be the standard basis for
$\mathbb{R}^n$, where $n$ is a fixed integer larger than $2$.
In this section, $D$ denotes a subset of $\mathbb{R}^n$ that
only satisfies $\{ 0, e_1, e_2, \ldots, e_n \} \subset D$,
whether bounded or not.

Lemma \ref{lem:trans} implies that there exists an orthogonal
matrix $\textbf{Q}$, with which we can express
$f(e_i) = (e_{i1}^\prime, e_{i2}^\prime, \ldots, e_{ii}^\prime,
           0, \ldots, 0)^{tr}$%, for all $i \in \{ 1, 2, \ldots, n \}$,
with respect to the new basis
$\{ \textbf{Q} e_1, \textbf{Q} e_2, \ldots, \textbf{Q} e_n \}$
instead of the standard basis $\{ e_1, e_2, \ldots, e_n \}$
and $e_{ii}^\prime \geq 0$ for all
$i \in \{ 1, 2, \ldots, n \}$.
For example, we can choose the orthogonal matrix $\textbf{Q}$
given in the proof of Lemma \ref{lem:trans} for this purpose.
Therefore, from now on, we will assume that
$f(e_i) = (e_{i1}^\prime, e_{i2}^\prime, \ldots, e_{ii}^\prime,
           0, \ldots, 0)$ with $e_{ii}^\prime \geq 0$, written
in row vector, without loss of generality.
%From now on, we will denote all vectors as row vectors for
%convenience.

In the statement of the following theorem, we define $\sigma$
as if we knew the values of $c_{(i+1)i}$ without knowing their
values in advance.
However, we note that in the proof of this theorem, we can
justify the definition of $\sigma$ by showing that the value
of each $c_{(i+1)i}$ is not greater than $9$
(ref. Remark \ref{rem:3.1} $(ii)$).
It is noted that the situation is similar for the $c_{ii}$'s.

%%%%%%%%%%%%%%%%%%%%%%%%%%%%%%%%%%%
%%%%%     Theorem 3.1     %%%%%%%%%
%%%%%%%%%%%%%%%%%%%%%%%%%%%%%%%%%%%

\begin{thm} \label{thm:ndlemma}
Let $\{ e_1, e_2, \ldots, e_n \}$ be the standard basis for
$\mathbb{R}^n$, where $\mathbb{R}^n$ denotes the $n$-dimensional
Euclidean space for a fixed integer $n \geq 3$, let $D$ be a
subset of $\mathbb{R}^n$ satisfying
$\{ 0, e_1, e_2, \ldots, e_n \} \subset D$, and let
$f : D \to \mathbb{R}^n$ be a function that satisfies
$f(0) = 0$ and
\begin{eqnarray} \label{eq:2021-1-13}
\big| \| f(x) - f(y) \| - \| x - y \| \big| \leq \varepsilon
\end{eqnarray}
for all $x, y \in \{ 0, e_1, e_2, \ldots, e_n \}$ and for some
constant $\varepsilon$ with
$0 < \varepsilon <
 \min \big\{ \frac{1}{\sigma},\,
             \min\limits_{1 \leq i \leq n} \frac{1}{2c_{ii}},\,
             \frac{1}{12}
      \big\}$,
where $\sigma$ is defined as
$\sigma = \sum\limits_{i=1}^{n-1} c_{(i+1)i}^2$ and $c_{ii}$
and $c_{(i+1)i}$ will be determined by the formulas
$(\ref{eq:diffeq1p})$ and $(\ref{eq:diffeq2p})$, respectively.
Then there exist positive integers $c_{ij}$,
$i, j \in \{ 1, 2, \ldots, n \}$ with $j \leq i$, such that
\begin{eqnarray} \label{eq:2021-1-22}
\left\{
\begin{array}{cl}
   -c_{ij} \varepsilon \leq e_{ij}^\prime
   \leq c_{ij} \varepsilon
   & (\mbox{for~$i > j$}), \vspace*{2mm}\\
   1 - c_{ii} \varepsilon \leq e_{ii}^\prime
   \leq 1 + \varepsilon
   & (\mbox{for~$i = j$})
\end{array}
\right.
\end{eqnarray}
and such that the $c_{ij}$ satisfy the equations in
$(\ref{eq:2021-1-23})$ for all $i, j \in \{ 1, 2, \ldots, n \}$
with $j \leq i$.
%In particular, each $c_{ij}$ is independent of $\varepsilon$.
\end{thm}
\vspace*{3mm}

\noindent
{\em Proof}.
$(a)$ Using inequality (\ref{eq:2021-1-13}) and by assumption
$f(0) = 0$, we have
\begin{eqnarray*}
\big| \| f(e_j) \| - 1 \big| \leq \varepsilon
~~~\mbox{and}~~~
\big| \| f(e_k) - f(e_\ell) \| - \sqrt{2} \big|
\leq \varepsilon
\end{eqnarray*}
for any $j, k, \ell \in \{ 1, 2, \ldots, n \}$ with
$k < \ell$.
Since
$f(e_j) = (e_{j1}^\prime, \ldots, e_{jj}^\prime, 0, \ldots, 0)$
for all $j \in \{ 1, 2, \ldots, n \}$ and $\| \cdot \|$ is the
Euclidean norm on $\mathbb{R}^n$, from the last inequalities
we get the following two inequalities, which are equivalent to
the inequality (\ref{eq:2021-1-13}) for
$x, y \in \{ 0, e_1, e_2, \ldots, e_n \}$:
\begin{eqnarray} \label{eq:el1}
(1 - \varepsilon)^2
\leq \sum_{i=1}^j e_{ji}^{\prime 2}
\leq (1 + \varepsilon)^2
\end{eqnarray}
for each $j \in \{ 1, 2, \ldots, n \}$ and
\begin{eqnarray} \label{eq:elk}
\big( \sqrt{2} - \varepsilon \big)^2
\leq \sum_{i=1}^k e_{ki}^{\prime 2} -
     \sum_{i=1}^k 2 e_{ki}^\prime e_{\ell i}^\prime +
     \sum_{i=1}^\ell e_{\ell i}^{\prime 2}
\leq \big( \sqrt{2} + \varepsilon \big)^2
\end{eqnarray}
for every $k, \ell \in \{ 1, 2, \ldots, n \}$ with $k < \ell$.
From now on, we will prove this theorem by using inequalities
(\ref{eq:el1}) and (\ref{eq:elk}) instead of inequality
(\ref{eq:2021-1-13}).

$(b)$ Now we apply the `main' induction on $m$ to prove the
array of equations presented in (\ref{eq:2021-1-23}).
Proving the array of (\ref{eq:2021-1-23}) is the most important
and longest part of this proof.

$(b.1)$ According to Lemma \ref{lem:trans}, $e_{11}^\prime$ is
a nonnegative real number, so setting $j = 1$ in (\ref{eq:el1})
gives us the inequality,
$1 - \varepsilon \leq e_{11}^\prime \leq 1 + \varepsilon$, and
we select $c_{11} = 1$ as the smallest positive integer that
satisfies the following inequality
\begin{eqnarray} \label{eq:2021-3-31}
1 - c_{11} \varepsilon \leq 1 - \varepsilon \leq e_{11}^\prime
\leq 1 + \varepsilon.
\end{eqnarray}
This fact guarantees the existence of $c_{11}$ satisfying the
second condition of (\ref{eq:2021-1-22}) for $i = j = 1$.
(We note that $c_{11}$ must not be necessarily the smallest
positive integer satisfying the inequality (\ref{eq:2021-3-31}),
for example, $c_{11} = 2$ is not wrong, but if possible, the
smaller the $c_{11}$ is, the better.)
If we set $j = 2$ in (\ref{eq:el1}) and put $k = 1$ and
$\ell = 2$ in (\ref{eq:elk}) and then combine the resulting
inequalities, then we get
\begin{eqnarray*}
\frac{-(2c_{11} + 2 + 2 \sqrt{2}) \varepsilon +
      c_{11}^2 \varepsilon^2}
     {2 (1 - c_{11} \varepsilon)}
\leq e_{21}^\prime
\leq \frac{(4 + 2 \sqrt{2}) \varepsilon + \varepsilon^2}
          {2 (1 - c_{11} \varepsilon)}
\end{eqnarray*}
and we can choose $c_{21} = 4$ as the smallest positive integer
satisfying the condition
\begin{eqnarray} \label{eq:2021-2-15a}
-c_{21} \varepsilon \leq
\frac{-(2c_{11} + 2 + 2 \sqrt{2}) \varepsilon +
      c_{11}^2 \varepsilon^2}
     {2 (1 - c_{11} \varepsilon)}
\leq e_{21}^\prime
\leq \frac{(4 + 2 \sqrt{2}) \varepsilon + \varepsilon^2}
          {2 (1 - c_{11} \varepsilon)}
\leq c_{21} \varepsilon.
\end{eqnarray}
Obviously, this fact confirms the existence of $c_{21}$
satisfying the first condition of (\ref{eq:2021-1-22}) for
$i = 2$ and $j = 1$.
(We note that $c_{ij}$ must not be necessarily the smallest
positive integer under given conditions, for example,
$c_{21} = 5$ is not bad, but if possible, the smaller the
$c_{ij}$ is, the better.
Indeed, assuming that the $c_{ij}$ are the smallest positive
integers that satisfy some given conditions, we can find the
unique $c_ {ij}$ (for $i > j$), which make it easier to prove
the array of equations presented in (\ref{eq:2021-1-23}).
This is one of the reasons we want $c_{ij}$ to be the smallest
positive integer.)

Furthermore, if we set $j = 3$ in (\ref{eq:el1}) and put
$k = 1$ and $\ell = 3$ in (\ref{eq:elk}) and then combine the
resulting inequalities, then we get the inner part of the
following inequalities
\begin{eqnarray} \label{eq:2021-2-22}
-c_{31} \varepsilon \leq
\frac{-(2c_{11} + 2 + 2 \sqrt{2}) \varepsilon +
      c_{11}^2 \varepsilon^2}
     {2 (1 - c_{11} \varepsilon)}
\leq e_{31}^\prime
\leq \frac{(4 + 2 \sqrt{2}) \varepsilon + \varepsilon^2}
          {2 (1 - c_{11} \varepsilon)}
\leq c_{31} \varepsilon
\end{eqnarray}
and we select $c_{31}$ as the smallest positive integer that
satisfies the outermost inequalities of (\ref{eq:2021-2-22}).
From this fact, we confirm the existence of $c_{31}$ that
satisfies the first condition in (\ref{eq:2021-1-22}).
Since $c_{21}$ and $c_{31}$ are assumed to be the smallest
positive integers satisfying the outermost inequalities of
(\ref{eq:2021-2-15a}) and (\ref{eq:2021-2-22}), respectively,
we come to the conclusion that $c_{31} = c_{21} = 4$.

Moreover, using (\ref{eq:el1}) with $j = 2$ and by a direct
calculation, since $0 < \varepsilon < \frac{1}{12}$, we get
\begin{eqnarray*}
1 - 2 \varepsilon - \frac{19}{12} \varepsilon + 4 \varepsilon^2
<    1 - 2 \varepsilon - 19 \varepsilon^2 + 4 \varepsilon^2
=    (1 - \varepsilon)^2 - c_{21}^2 \varepsilon^2
\leq e_{22}^{\prime 2}
\leq (1 + \varepsilon)^2
\end{eqnarray*}
and, if possible, we determine $c_{22}$ as the smallest
positive integer which satisfies the condition
\begin{eqnarray*}
(1 - c_{22} \varepsilon)^2
\leq 1 - 2 \varepsilon - \frac{19}{12} \varepsilon +
     4 \varepsilon^2
\leq e_{22}^{\prime 2}
\leq (1 + \varepsilon)^2.
\end{eqnarray*}
The last inequalities assures the existence of $c_{22} = 2$
satisfying the second condition of (\ref{eq:2021-1-22}) for
$i = j = 2$.

In a similar way, putting $k = 2$ and $\ell = 3$ in
(\ref{eq:elk}) and a routine calculation show the existence of
$c_{32}$ satisfying the first condition of (\ref{eq:2021-1-22}),
for example, $c_{32} = 6$.
Analogously, since $0 < \varepsilon < \frac{1}{12}$, the
inequality (\ref{eq:el1}) with $j = 3$ yields the inequality,
$1 - 2 \varepsilon - 51 \varepsilon^2 \leq e_{33}^{\prime 2}
 \leq (1 + \varepsilon)^2$, and we can choose the smallest
positive integer $c_{33}$ that satisfies
\begin{eqnarray*}
(1 - c_{33} \varepsilon)^2
\leq 1 - 2 \varepsilon - 51 \varepsilon^2 \leq e_{33}^{\prime 2}
\leq (1 + \varepsilon)^2,
\end{eqnarray*}
which shows that $c_{33}$ exists that satisfies the second
condition of (\ref{eq:2021-1-22}).
More directly, we can choose $c_{33} = 4$.
Therefore, all the integers $c_{ij}$ considered in $(b.1)$
satisfy the conditions in (\ref{eq:2021-1-22}) and
(\ref{eq:2021-1-23}) for $n = 3$.
By doing this, we start the induction (with $m = 3$).

$(b.2)$ Induction hypothesis.
Let $m$ be some integer satisfying $3 \leq m < n$.
It is assumed that the smallest positive integers $c_{ij}$,
$i, j \in \{ 1, 2, \ldots, m \}$ with $j \leq i$, were found
by the methods we did in the subsection $(b.1)$, and that these
positive integers satisfy the following inequalities
\begin{eqnarray*}
\left\{
\begin{array}{cl}
   -c_{ij} \varepsilon \leq e_{ij}^\prime
   \leq c_{ij} \varepsilon
   & (\mbox{for~$i > j$}), \vspace*{2mm}\\
   1 - c_{ii} \varepsilon \leq e_{ii}^\prime
   \leq 1 + \varepsilon
   & (\mbox{for~$i = j$})
\end{array}
\right.
\end{eqnarray*}
as well as the array of equations
\begin{eqnarray*} %\label{eq:indass}
\left\{
\begin{array}{ccccccccccccc}
c_{m1}                & \!\!\!\!\!=\!\!\!\!\! &
c_{(m-1)1}            & \!\!\!\!\!=\!\!\!\!\! &
c_{(m-2)1}            & \!\!\!\!\!=\!\!\!\!\! &
~~\ldots~~            & \!\!\!\!\!=\!\!\!\!\! &
c_{41}                & \!\!\!\!\!=\!\!\!\!\! &
c_{31}                & \!\!\!\!\!=\!\!\!\!   &
c_{21}, \\
c_{m2}                & \!\!\!\!\!=\!\!\!\!\! &
c_{(m-1)2}            & \!\!\!\!\!=\!\!\!\!\! &
c_{(m-2)2}            & \!\!\!\!\!=\!\!\!\!\! &
~~\ldots~~            & \!\!\!\!\!=\!\!\!\!\! &
c_{42}                & \!\!\!\!\!=\!\!\!\!   &
c_{32}, \\
c_{m3}                & \!\!\!\!\!=\!\!\!\!\! &
c_{(m-1)3}            & \!\!\!\!\!=\!\!\!\!\! &
c_{(m-2)3}            & \!\!\!\!\!=\!\!\!\!\! &
~~\ldots~~            & \!\!\!\!\!=\!\!\!\!   &
c_{43}, \\
\vdots                &                       &
\vdots                &                       &
\vdots \\
c_{m(m-3)}            & \!\!\!\!\!=\!\!\!\!\! &
c_{(m-1)(m-3)}        & \!\!\!\!\!=\!\!\!\!\! &
c_{(m-2)(m-3)}, \\
c_{m(m-2)}            & \!\!\!\!\!=\!\!\!\!\! &
c_{(m-1)(m-2)}, \\
c_{m(m-1)}
\end{array}
\right.
\end{eqnarray*}
The last line in the above array consisting of only $c_{m(m-1)}$
means that there exists the smallest possible positive integer
$c_{m(m-1)}$ that satisfies
$-c_{m(m-1)} \varepsilon \leq e_{m(m-1)}^\prime
 \leq c_{m(m-1)} \varepsilon$.

$(b.3)$ We let $j = m+1$ in (\ref{eq:el1}) and $\ell = m+1$ in
(\ref{eq:elk}) to get
\begin{eqnarray} \label{eq:en1}
(1 - \varepsilon)^2
\leq \sum_{i=1}^{m+1} e_{(m+1)i}^{\prime 2}
\leq (1 + \varepsilon)^2
\end{eqnarray}
and
\begin{eqnarray} \label{eq:enk}
\big( \sqrt{2} - \varepsilon \big)^2
\leq \sum_{i=1}^k e_{ki}^{\prime 2} -
     \sum_{i=1}^k 2 e_{ki}^\prime e_{(m+1)i}^\prime +
     \sum_{i=1}^{m+1} e_{(m+1)i}^{\prime 2}
\leq \big( \sqrt{2} + \varepsilon \big)^2
\end{eqnarray}
for every $k \in \{ 1, 2, \ldots, m \}$.

Similar to what we did to get (\ref{eq:2021-2-15a}), the
inequalities (\ref{eq:2021-3-31}), (\ref{eq:en1}), and
(\ref{eq:enk}) with $k = 1$ yield the inner ones of the
following inequalities
\begin{eqnarray} \label{eq:2021-2-15b}
-c_{(m+1)1} \varepsilon \leq
\frac{-(2c_{11} + 2 + 2 \sqrt{2}) \varepsilon +
      c_{11}^2 \varepsilon^2}
     {2 (1 - c_{11} \varepsilon)}
\leq e_{(m+1)1}^\prime
\leq \frac{(4 + 2 \sqrt{2}) \varepsilon + \varepsilon^2}
          {2 (1 - c_{11} \varepsilon)}
\leq c_{(m+1)1} \varepsilon,
\end{eqnarray}
and we find the smallest positive integer $c_{(m+1)1}$
satisfying the outermost inequalities of (\ref{eq:2021-2-15b}).
By comparing both inequalities (\ref{eq:2021-2-15a}) and
(\ref{eq:2021-2-15b}), we may conclude that
$c_{(m+1)1} = c_{21}$, with which we initiate an `inner'
induction that is subordinate to the main induction.
%%%%%%%%%%%%%%%%%%%%%%%%%%%%%%%%%%%%%%%%%%%%%%%%%%%%%%%%%%%
%%%   We just began small inner induction     %%%%%%%%%%%%%
%%%%%%%%%%%%%%%%%%%%%%%%%%%%%%%%%%%%%%%%%%%%%%%%%%%%%%%%%%%

$(b.3.1)$ We choose some $k \in \{ 2, 3, \ldots, m \}$ and
assume that
$-c_{(m+1)i} \varepsilon \leq e_{(m+1)i}^\prime \leq c_{(m+1)i}
 \varepsilon$ and $c_{(m+1)i} = c_{(i+1)i}$ for each
$i \in \{ 1, 2, \ldots, k-1 \}$.
This is the hypothesis for our inner induction on $i$ that
operates inside the main induction on $m$.
Based on this hypothesis, we will prove that there exists a
positive integer $c_{(m+1)k}$ that satisfies
$-c_{(m+1)k} \varepsilon \leq e_{(m+1)k}^\prime \leq c_{(m+1)k}
 \varepsilon$ as well as $c_{(m+1)k} = c_{(k+1)k}$.
Roughly speaking, this inner induction works to expand each
row of (\ref{eq:2021-1-23}) horizontally.

$(b.3.2)$ It follows from (\ref{eq:enk}) that
\begin{eqnarray} \label{eq:2021-1-15}
& \big( \sqrt{2} - \varepsilon \big)^2 -
  {\displaystyle
  \sum_{i=1}^k e_{ki}^{\prime 2} +
  \sum_{i=1}^{k-1} 2 e_{ki}^\prime e_{(m+1)i}^\prime -
  \sum_{i=1}^{m+1} e_{(m+1)i}^{\prime 2}}
& \nonumber \vspace*{2mm}\\
& \leq -2 e_{kk}^\prime e_{(m+1)k}^\prime \leq
& \vspace*{2mm}\\
& \big( \sqrt{2} + \varepsilon \big)^2 -
  {\displaystyle
  \sum_{i=1}^k e_{ki}^{\prime 2} +
  \sum_{i=1}^{k-1} 2 e_{ki}^\prime e_{(m+1)i}^\prime -
  \sum_{i=1}^{m+1} e_{(m+1)i}^{\prime 2}}
& \nonumber
\end{eqnarray}
for any $k \in \{ 2, 3, \ldots, m \}$.
On the other hand, by (\ref{eq:el1}) and (\ref{eq:en1}), we
have
\begin{eqnarray*}
(1 - \varepsilon)^2
\leq \sum_{i=1}^k e_{ki}^{\prime 2}
\leq (1 + \varepsilon)^2
~~~\mbox{and}~~~
(1 - \varepsilon)^2
\leq \sum_{i=1}^{m+1} e_{(m+1)i}^{\prime 2}
\leq (1 + \varepsilon)^2
\end{eqnarray*}
for each $k \in \{ 2, 3, \ldots, m \}$.
Moreover, it follows from the hypotheses $(b.2)$ and $(b.3.1)$
that
\begin{eqnarray*}
-\sum_{i=1}^{k-1} 2 c_{ki} c_{(i+1)i} \varepsilon^2
=    -\sum_{i=1}^{k-1} 2 c_{ki} c_{(m+1)i} \varepsilon^2
\leq \sum_{i=1}^{k-1} 2 e_{ki}^\prime e_{(m+1)i}^\prime
\leq \sum_{i=1}^{k-1} 2 c_{ki} c_{(m+1)i} \varepsilon^2
=    \sum_{i=1}^{k-1} 2 c_{ki} c_{(i+1)i} \varepsilon^2
\end{eqnarray*}
for all $k \in \{ 2, 3, \ldots, m \}$.

Since $c_{kk} \varepsilon < \frac{1}{2}$ and $e_{kk}^\prime > 0$
by $(b.2)$, we use (\ref{eq:2021-1-15}) and the last
inequalities to get the inner ones of the following inequalities
\begin{eqnarray} \label{eq:2021-1-16a}
& -c_{(m+1)k} \varepsilon \leq
  {\displaystyle
  \frac{1}{2 e_{kk}^\prime}
  \left( -\big( 4 + 2 \sqrt{2} \big) \varepsilon +
         \varepsilon^2 -
         2 \sum_{i=1}^{k-1} c_{ki} c_{(i+1)i} \varepsilon^2
  \right)}
& \nonumber \vspace*{2mm}\\
& \leq e_{(m+1)k}^\prime \leq
& \vspace*{2mm}\\
& {\displaystyle
  \frac{1}{2 e_{kk}^\prime}
  \left( \big( 4 + 2 \sqrt{2} \big) \varepsilon +
         \varepsilon^2 +
         {\displaystyle
         2 \sum_{i=1}^{k-1} c_{ki} c_{(i+1)i} \varepsilon^2}
  \right)}
  \leq c_{(m+1)k} \varepsilon
& \nonumber
\end{eqnarray}
for all $k \in \{ 2, 3, \ldots, m \}$ and we can select the
smallest positive integer $c_{(m+1)k}$ that satisfies the
outermost inequalities of (\ref{eq:2021-1-16a}).
Similarly, by (\ref{eq:el1}) and (\ref{eq:elk}) with
$\ell = k + 1$, a routine calculation yields
\begin{eqnarray} \label{eq:2021-2-15c}
& -c_{(k+1)k} \varepsilon \leq
  {\displaystyle
  \frac{1}{2 e_{kk}^\prime}
  \left( -\big( 4 + 2 \sqrt{2} \big) \varepsilon +
         \varepsilon^2 -
         2 \sum_{i=1}^{k-1} c_{ki} c_{(k+1)i} \varepsilon^2
  \right)}
& \nonumber \vspace*{2mm}\\
& \leq e_{(k+1)k}^\prime \leq
& \vspace*{2mm}\\
& {\displaystyle
  \frac{1}{2 e_{kk}^\prime}
  \left( \big( 4 + 2 \sqrt{2} \big) \varepsilon +
         \varepsilon^2 +
         {\displaystyle
         2 \sum_{i=1}^{k-1} c_{ki} c_{(k+1)i} \varepsilon^2}
  \right)}
  \leq c_{(k+1)k} \varepsilon,
& \nonumber
\end{eqnarray}
where $c_{(k+1)k}$ is the smallest positive integer that
satisfies the outmost conditions of (\ref{eq:2021-2-15c}).
We note by $(b.2)$ and $(b.3.1)$ that $c_{(k+1)i} = c_{(i+1)i}$
for every integer $i$ satisfying $0 < i < k$.
Comparing (\ref{eq:2021-1-16a}) and (\ref{eq:2021-2-15c}), we
conclude that $c_{(m+1)k} = c_{(k+1)k}$ for each
$k \in \{ 2, 3, \ldots, m \}$.
Furthermore, referring to the subsection $(b.3)$, we see that
$c_{(m+1)k} = c_{(k+1)k}$ holds for all
$k \in \{ 1, 2, \ldots, m \}$, which proves the truth of the
first column of the array of equations in subsection $(b.3.3)$
below.

Moreover, the inequality (\ref{eq:el1}) with $j = m+1$ yields
\begin{eqnarray} \label{eq:2021-2-20}
(1 - \varepsilon)^2 - \sum_{i=1}^m e_{(m+1)i}^{\prime 2}
\leq e_{(m+1)(m+1)}^{\prime 2}
\leq (1 + \varepsilon)^2 - \sum_{i=1}^m e_{(m+1)i}^{\prime 2}.
\end{eqnarray}
Since $0 < \varepsilon < \frac{1}{\sigma}$, $m < n$, and
$0 < \varepsilon < \frac{1}{12}$, it follows from
(\ref{eq:2021-2-20}) and some manipulation that
\begin{eqnarray*}
\begin{split}
(1 + \varepsilon)^2
& \geq e_{(m+1)(m+1)}^{\prime 2}
  \geq (1 - \varepsilon)^2 -
       \sum_{i=1}^m c_{(m+1)i}^2 \varepsilon^2
  =    (1 - \varepsilon)^2 -
       \sum_{i=1}^m c_{(i+1)i}^2 \varepsilon^2 \\
& \geq (1 - \varepsilon)^2 - \sigma \varepsilon^2
  >    1 - 3 \varepsilon + \varepsilon^2
  >    9 \varepsilon + \varepsilon^2
%  >    109 \varepsilon^2 \\
  >    100 \varepsilon^2.
\end{split}
\end{eqnarray*}
We see that $0 < (1 - c \varepsilon)^2 \leq 100 \varepsilon^2$
whenever $c$ is a positive integer satisfying
$\frac{1}{\varepsilon} - 10 \leq c < \frac{1}{\varepsilon}$.
This fact shows the existence of $c_{(m+1)(m+1)}$ that
satisfies the second condition of (\ref{eq:2021-1-22}).

$(b.3.3)$ We just proved in the subsections from $(b.2)$ to
$(b.3.2)$ that there exist positive integers $c_{ij}$,
$i, j \in \{ 1, 2, \ldots, m+1 \}$ with $j \leq i$, such that
\begin{eqnarray*}
\left\{
\begin{array}{cl}
   -c_{ij} \varepsilon \leq e_{ij}^\prime
   \leq c_{ij} \varepsilon
   & (\mbox{for~$i > j$}), \vspace*{2mm}\\
   1 - c_{ii} \varepsilon \leq e_{ii}^\prime
   \leq 1 + \varepsilon
   & (\mbox{for~$i = j$}),
\end{array}
\right.
\end{eqnarray*}
and the $c_{ij}$'s satisfy
\begin{eqnarray*}
\left\{
\begin{array}{ccccccccccccc}
c_{(m+1)1}            & \!\!\!\!\!=\!\!\!\!\! &
c_{m1}                & \!\!\!\!\!=\!\!\!\!\! &
c_{(m-1)1}            & \!\!\!\!\!=\!\!\!\!\! &
~~\ldots~~            & \!\!\!\!\!=\!\!\!\!\! &
c_{41}                & \!\!\!\!\!=\!\!\!\!\! &
c_{31}                & \!\!\!\!\!=\!\!\!\!   &
c_{21}, \\
c_{(m+1)2}            & \!\!\!\!\!=\!\!\!\!\! &
c_{m2}                & \!\!\!\!\!=\!\!\!\!\! &
c_{(m-1)2}            & \!\!\!\!\!=\!\!\!\!\! &
~~\ldots~~            & \!\!\!\!\!=\!\!\!\!\! &
c_{42}                & \!\!\!\!\!=\!\!\!\!   &
c_{32}, \\
c_{(m+1)3}            & \!\!\!\!\!=\!\!\!\!\! &
c_{m3}                & \!\!\!\!\!=\!\!\!\!\! &
c_{(m-1)3}            & \!\!\!\!\!=\!\!\!\!\! &
~~\ldots~~            & \!\!\!\!\!=\!\!\!\!   &
c_{43}, \\
\vdots                &                       &
\vdots                &                       &
\vdots \\
c_{(m+1)(m-2)}        & \!\!\!\!\!=\!\!\!\!\! &
c_{m(m-2)}            & \!\!\!\!\!=\!\!\!\!\! &
c_{(m-1)(m-2)}, \\
c_{(m+1)(m-1)}        & \!\!\!\!\!=\!\!\!\!\! &
c_{m(m-1)}, \\
c_{(m+1)m}
\end{array}
\right.
\end{eqnarray*}
The last row in the above array consisting of only $c_{(m+1)m}$
means that there exists the smallest possible positive integer
$c_{(m+1)m}$ that satisfies
$-c_{(m+1)m} \varepsilon \leq e_{(m+1)m}^\prime
 \leq c_{(m+1)m} \varepsilon$.
(We can check in $(b.3.2)$ the truth of equations in the first
column of the above array.
Moreover, we remember that we have assumed in $(b.2)$ that the
rest of equations in the array are true.)

$(b.4)$ Altogether, by the main induction conclusion on $m$
($3 \leq m < n$), we may conclude that there exist positive
integers $c_{ij}$, $i, j \in \{ 1, 2, \ldots, n \}$ with
$j \leq i$, such that each inequality in (\ref{eq:2021-1-22})
holds true and the $c_{ij}$'s satisfy
\begin{eqnarray} \label{eq:2021-1-23}
\left\{
\begin{array}{ccccccccccccc}
c_{n1}                & \!\!\!\!\!=\!\!\!\!\! &
c_{(n-1)1}            & \!\!\!\!\!=\!\!\!\!\! &
c_{(n-2)1}            & \!\!\!\!\!=\!\!\!\!\! &
~~\ldots~~            & \!\!\!\!\!=\!\!\!\!\! &
c_{41}                & \!\!\!\!\!=\!\!\!\!\! &
c_{31}                & \!\!\!\!\!=\!\!\!\!   &
c_{21}, \\
c_{n2}                & \!\!\!\!\!=\!\!\!\!\! &
c_{(n-1)2}            & \!\!\!\!\!=\!\!\!\!\! &
c_{(n-2)2}            & \!\!\!\!\!=\!\!\!\!\! &
~~\ldots~~            & \!\!\!\!\!=\!\!\!\!\! &
c_{42}                & \!\!\!\!\!=\!\!\!\!   &
c_{32}, \\
c_{n3}                & \!\!\!\!\!=\!\!\!\!\! &
c_{(n-1)3}            & \!\!\!\!\!=\!\!\!\!\! &
c_{(n-2)3}            & \!\!\!\!\!=\!\!\!\!\! &
~~\ldots~~            & \!\!\!\!\!=\!\!\!\!   &
c_{43}, \\
\vdots                &                       &
\vdots                &                       &
\vdots \\
c_{n(n-3)}            & \!\!\!\!\!=\!\!\!\!\! &
c_{(n-1)(n-3)}        & \!\!\!\!\!=\!\!\!\!\! &
c_{(n-2)(n-3)}, \\
c_{n(n-2)}            & \!\!\!\!\!=\!\!\!\!\! &
c_{(n-1)(n-2)}, \\
c_{n(n-1)}
\end{array}
\right.
\end{eqnarray}
which completes the first part of our proof.
We remark that the last row `$c_{n(n-1)}$' in the above array
implies that there is an integer $c_{n(n-1)} > 0$ satisfying
$-c_{n(n-1)} \varepsilon \leq e_{n(n-1)}^\prime
 \leq c_{n(n-1)} \varepsilon$.

$(c)$ Now we will introduce efficient methods to estimate the
positive integers $c_{jj}$ and $c_{k(k-1)}$ for every
$j \in \{ 1, 2, \ldots, n \}$ and $k \in \{ 2, 3, \ldots, n \}$.

$(c.1)$ We note that the inequality (\ref{eq:el1}) holds true
for all $j \in \{ 1, 2, \ldots, n \}$.
Since
$-c_{ji} \varepsilon \leq e_{ji}^\prime \leq c_{ji} \varepsilon$
for any $i, j \in \{ 1, 2, \ldots, n \}$ with $i < j$, we
determine the $c_{jj}$ as the smallest possible positive
integer that satisfies
\begin{eqnarray} \label{eq:2021-2-28}
(1 - c_{jj} \varepsilon)^2
\leq (1 - \varepsilon)^2 - \sum_{i=1}^{j-1} c_{ji}^2 \varepsilon^2
\leq (1 - \varepsilon)^2 - \sum_{i=1}^{j-1} e_{ji}^{\prime 2}
\leq e_{jj}^{\prime 2}
\leq (1 + \varepsilon)^2
\end{eqnarray}
for all $j \in \{ 2, 3, \ldots, n \}$.
However, since the previous inequality is inefficient in
practical calculations, we introduce a more practical inequality
even at the expense of the smallest property of $c_{jj}$.
Instead of (\ref{eq:2021-2-28}), we will determine the $c_{jj}$
as the smallest positive integer that satisfies the new
condition
\begin{eqnarray} \label{eq:diffeq1p}
(3c_{jj} - 1) (c_{jj} - 1) \geq \sum_{i=1}^{j-1} c_{(i+1)i}^2
\end{eqnarray}
for all $j \in \{ 2, 3, \ldots, n \}$, which proves the
existence of $c_{jj}$.
Indeed, since $0 < \varepsilon < \frac{1}{2c_{jj}}$, we have
\begin{eqnarray*}
(1 - c_{jj} \varepsilon)^2
= 1 - 2 \varepsilon c_{jj} + \varepsilon^2 c_{jj}^2
< 1 - 2 \varepsilon c_{jj} + \frac{1}{2} \varepsilon c_{jj}
= 1 - \frac{3}{2} \varepsilon c_{jj}.
\end{eqnarray*}
It further follows from (\ref{eq:2021-1-23}),
(\ref{eq:diffeq1p}), and the last inequality that
\begin{eqnarray*}
\begin{split}
(1 - c_{jj} \varepsilon)^2
& <    1 - \frac{3}{2} \varepsilon c_{jj}
  =    1 - \frac{\varepsilon}{2 c_{jj}} \cdot 3 c_{jj}^2
  \leq 1 - \frac{\varepsilon}{2c_{jj}}
           \left( \sum_{i=1}^{j-1} c_{(i+1)i}^2 - 1 + 4 c_{jj}
           \right) \\
& =    1 - 2 \varepsilon - \frac{\varepsilon}{2c_{jj}}
       \left( \sum_{i=1}^{j-1} c_{(i+1)i}^2 - 1 \right)
  \leq 1 - 2 \varepsilon - \varepsilon^2
       \left( \sum_{i=1}^{j-1} c_{(i+1)i}^2 - 1 \right) \\
& =    (1 - \varepsilon)^2 -
       \sum_{i=1}^{j-1} c_{(i+1)i}^2 \varepsilon^2
  =    (1 - \varepsilon)^2 -
       \sum_{i=1}^{j-1} c_{ji}^2 \varepsilon^2,
\end{split}
\end{eqnarray*}
which implies the validity of (\ref{eq:2021-2-28}).

$(c.2)$ We note that the inequality (\ref{eq:elk}) holds for
all $k, \ell \in \{ 1, 2, \ldots, n \}$ with $k < \ell$.
If we set $\ell = k + 1$ in (\ref{eq:elk}) and make some
manipulations, then we obtain the inner ones of the following
inequalities
\begin{eqnarray} \label{eq:2021-1-16c}
& -c_{(k+1)k} \varepsilon \leq
  {\displaystyle
  \frac{1}{2 e_{kk}^\prime}
  \left( -\big( 4 + 2 \sqrt{2} \big) \varepsilon +
         \varepsilon^2 -
         2 \sum_{i=1}^{k-1} c_{ki} c_{(k+1)i} \varepsilon^2
  \right)}
& \nonumber \vspace*{2mm}\\
& \leq e_{(k+1)k}^\prime \leq
& \vspace*{2mm}\\
& {\displaystyle
  \frac{1}{2 e_{kk}^\prime}
  \left( \big( 4 + 2 \sqrt{2} \big) \varepsilon +
         \varepsilon^2 +
         {\displaystyle
         2 \sum_{i=1}^{k-1} c_{ki} c_{(k+1)i} \varepsilon^2}
  \right)}
  \leq c_{(k+1)k} \varepsilon
& \nonumber
\end{eqnarray}
for every $k \in \{ 1, 2, \ldots, n-1 \}$.
And then, we choose the smallest positive integer $c_{(k+1)k}$
that satisfies the outermost inequalities of
(\ref{eq:2021-1-16c}).

The inequalities (\ref{eq:2021-1-22}) and (\ref{eq:2021-1-16c})
show the existence of the smallest positive integer $c_{(k+1)k}$
that satisfies
\begin{eqnarray} \label{eq:diffeq2p}
\frac{1}{2 (1 - c_{kk} \varepsilon)}
\left( 4 + 2 \sqrt{2} + \varepsilon +
       2 \sum_{i=1}^{k-1} c_{ki} c_{(k+1)i} \varepsilon
\right)
\leq c_{(k+1)k}
\end{eqnarray}
for $k \in \{ 1, 2, \ldots, n-1 \}$:
Since $0 < \varepsilon < \frac{1}{2c_{kk}}$, we see that
$\frac{1}{2 (1 - c_{kk} \varepsilon)} < 1$.
Further, since $c_{ki} = c_{(k+1)i} = c_{(i+1)i}$ for any
$i \in \{ 1, 2, \ldots, k-1 \}$ and
$0 < \varepsilon < \frac{1}{\sigma}$, we know that
$\sum\limits_{i=1}^{k-1} c_{ki} c_{(k+1)i} \varepsilon \leq
 \sigma \varepsilon < 1$.
Thus, we get
\begin{eqnarray*}
\frac{1}{2 (1 - c_{kk} \varepsilon)}
\left( 4 + 2 \sqrt{2} + \varepsilon +
       2 \sum_{i=1}^{k-1} c_{ki} c_{(k+1)i} \varepsilon
\right)
< 6 + 2 \sqrt{2} + \varepsilon
< 9,
\end{eqnarray*}
which together with (\ref{eq:diffeq2p}) assures the existence
of $c_{(k+1)k}$ with $0 < c_{(k+1)k} \leq 9$ for all
$k \in \{ 1, 2, \ldots, n-1 \}$.
\hspace*{\fill}$\Box$
\vspace*{3mm}
%%%%%%%%%%%%%%%%%%%%%%%%%%%%%%%%%%
%%%   End of Theorem 3.1   %%%%%%%
%%%%%%%%%%%%%%%%%%%%%%%%%%%%%%%%%%

%%%%%%%%%%%%%%%%%%%%%%%%%%%%%%%%%%%
%%%%%     Remark 3.1     %%%%%%%%%%
%%%%%%%%%%%%%%%%%%%%%%%%%%%%%%%%%%%

\begin{rem} \label{rem:3.1}
\begin{itemize}
\item[$(i)$]   The inequality $(\ref{eq:2021-1-13})$ is a
sufficient condition for inequalities in $(\ref{eq:2021-1-22})$,
and the inequalities in $(\ref{eq:2021-1-22})$ are necessary
conditions for the inequality $(\ref{eq:2021-1-13})$.

\item[$(ii)$]  In view of the last part of the proof of Theorem
\ref{thm:ndlemma}, we see that $0 < c_{(i+1)i} \leq 9$ for any
$i \in \{ 1, 2, \ldots, n-1 \}$.

\item[$(iii)$] We solve the quadratic inequality
$(\ref{eq:diffeq1p})$ with respect to $c_{jj}$ as follows:
\begin{eqnarray*}
c_{jj} \geq
\frac{1}{3}
\left( 2 + \sqrt{1 + 3 \sum\limits_{i=1}^{j-1} c_{(i+1)i}^2}\,
\right)
\end{eqnarray*}
for all $j \in \{ 1, 2, \ldots, n \}$.
\end{itemize}
\end{rem}

%%%%%%%%%%%%%%%%%%%%%%%%%%%%%%%%%%%%%%%%%%%%%%%%%%%%%%%%%%
%%%%%     4. Hyers-Ulam stability of isometries     %%%%%%
%%%%%%%%%%%%%%%%%%%%%%%%%%%%%%%%%%%%%%%%%%%%%%%%%%%%%%%%%%

\section{Hyers-Ulam stability of isometries on bounded domains}

The following theorem significantly improves the Fickett's
theorem by demonstrating the Hyers-Ulam stability of isometries
on the bounded domains.

As before, let $\{ e_1, e_2, \ldots, e_n \}$ be the standard
basis for $\mathbb{R}^n$.
Based on Lemma \ref{lem:trans}, we can assume that
$f(e_i) = (e_{i1}^\prime, e_{i2}^\prime, \ldots, e_{ii}^\prime,
           0, \ldots, 0)$ is written in row vector, where
$e_{ii}^\prime \geq 0$ for each $i \in \{ 1, 2, \ldots, n \}$.
We denote by $B_d(0)$ the closed ball of radius $d$ and
centered at the origin of $\mathbb{R}^n$, \textit{i.e.},
$B_d(0) = \{ x \in \mathbb{R}^n : \| x \| \leq d \}$.

%%%%%%%%%%%%%%%%%%%%%%%%%%%%%%%%%%%%%
%%%%%     Theorem 4.1     %%%%%%%%%%%
%%%%%%%%%%%%%%%%%%%%%%%%%%%%%%%%%%%%%

\begin{thm} \label{thm:ndtheorem}
Assume that $\mathbb{R}^n$ is the $n$-dimensional Euclidean
space for a given integer $n \geq 3$.
Let $D$ be a subset of $\mathbb{R}^n$ satisfying
$\{ 0, e_1, e_2, \ldots, e_n \} \subset D \subset B_d(0)$ for
some $d \geq 1$ and let $f : D \to \mathbb{R}^n$ be a function
satisfying $f(0) = 0$ and the inequality $(\ref{eq:2021-1-13})$
for all $x, y \in D$ and some constant $\varepsilon$ with
$0 < \varepsilon <
 \min \big\{ \frac{1}{\sigma},\,
             \min\limits_{1 \leq i \leq n} \frac{1}{2c_{ii}},\,
             \frac{1}{12}
     \big\}$,
where $\sigma = \sum\limits_{i=1}^{n-1} c_{(i+1)i}^2$ and the
$c_{ij}$, for all $i, j \in \{ 1, 2, \ldots, n \}$ with
$j \leq i$, are the positive integers estimated in Theorem
\ref{thm:ndlemma}.
Then there exists an isometry $U : D \to \mathbb{R}^n$ such
that
\begin{eqnarray*}
\| f(x) - U(x) \|
\leq \left[ \sum_{i=1}^n
            \left( \left( 2 + \sum_{j=1}^i c_{ij} \right) d +
                   4 + \sum_{j=1}^i c_{ij}
            \right)^2 \,
     \right]^{1/2} \varepsilon
\end{eqnarray*}
for all $x \in D$.
\end{thm}
\vspace*{3mm}

\noindent
{\em Proof}.
Let $\{ e_1, e_2, \ldots, e_n \}$ be the standard basis for
$\mathbb{R}^n$.
Based on Lemma \ref{lem:trans}, it can be assumed that
$f(e_i) = (e_{i1}^\prime, e_{i2}^\prime, \ldots, e_{ii}^\prime,
           0, \ldots, 0)$, where $e_{ii}^\prime \geq 0$ for
each $i \in \{ 1, 2, \ldots, n \}$.
For an arbitrary point $x = (x_1, x_2, \ldots, x_n)$ of $D$, let
$f(x) = (x_1^\prime, x_2^\prime, \ldots, x_n^\prime)$.
It follows from (\ref{eq:2021-1-13}) that
\begin{eqnarray*}
\big| \| f(x) \| - \| x \| \big| \leq \varepsilon
~~~\mbox{and}~~~
\big| \| f(x) - f(e_j) \| - \| x - e_j \| \big| \leq \varepsilon
\end{eqnarray*}
and hence, we have
\begin{eqnarray}
& \left| \left( {\displaystyle \sum_{i=1}^n x_i^{\prime 2}}
         \right)^{1/2} -
         \left( {\displaystyle \sum_{i=1}^n x_i^2} \right)^{1/2}
  \right| \leq \varepsilon, & \label{eq:2021-1-5an}
  \vspace*{3mm}\\
& \left| \left( {\displaystyle
                \sum_{i=1}^j (x_i^\prime - e_{ji}^\prime)^2 +
                \sum_{i=j+1}^n x_i^{\prime 2}}
         \right)^{1/2} -
         \left( {\displaystyle \sum_{i=1}^n x_i^2 -
                               2x_j + 1} \right)^{1/2}
  \right| \leq \varepsilon & \label{eq:2021-1-5dn}
\end{eqnarray}
for all $j \in \{ 1, 2, \ldots, n \}$.

It follows from (\ref{eq:2021-1-5an}) that
\begin{eqnarray} \label{eq:2021-1-6an}
\begin{split}
\left| \sum_{i=1}^n x_i^{\prime 2} - \sum_{i=1}^n x_i^2
\right|
& =    \left| \left( \sum_{i=1}^n x_i^{\prime 2} \right)^{1/2} +
              \left( \sum_{i=1}^n x_i^2 \right)^{1/2}
       \right|
       \left| \left( \sum_{i=1}^n x_i^{\prime 2} \right)^{1/2} -
              \left( \sum_{i=1}^n x_i^2 \right)^{1/2}
       \right| \\
& \leq (2d + 1) \varepsilon
\end{split}
\end{eqnarray}
since
$\sqrt{x_1^{\prime 2} + \cdots + x_n^{\prime 2}}
 \leq d + \varepsilon$,
$\sqrt{x_1^2 + \cdots + x_n^2} \leq d$, and
$0 < \varepsilon < 1$.
Similarly, it follows from (\ref{eq:2021-1-5dn}) that
\begin{eqnarray} \label{eq:2021-1-6bn}
\begin{split}
\left| \left( {\displaystyle
              \sum_{i=1}^j (x_i^\prime - e_{ji}^\prime)^2 +
              \sum_{i=j+1}^n x_i^{\prime 2}}
       \right) -
       \left( {\displaystyle \sum_{i=1}^n x_i^2 - 2x_j + 1}
       \right)
\right| \leq (2d + 3) \varepsilon
\end{split}
\end{eqnarray}
for all $j \in \{ 1, 2, \ldots, n \}$, since
$0 < \varepsilon < 1$ and
\begin{eqnarray*}
\left( \sum_{i=1}^j (x_i^\prime - e_{ji}^\prime)^2 +
       \sum_{i=j+1}^n x_i^{\prime 2}
\right)^{1/2} \leq d + 1 + \varepsilon
~~~\mbox{and}~~~
\left( \sum_{i=1}^n x_i^2 - 2x_j + 1 \right)^{1/2}
\leq d + 1.
\end{eqnarray*}

We use (\ref{eq:2021-1-6bn}) to get
\begin{eqnarray} \label{eq:2021-1-21}
& -(2d + 3) \varepsilon -
  {\displaystyle
  \left( \sum_{i=1}^n x_i^{\prime 2} - \sum_{i=1}^n x_i^2
  \right) + \sum_{i=1}^{j-1} 2 e_{ji}^\prime x_i^\prime +
  1 - \sum_{i=1}^j e_{ji}^{\prime 2}}
& \nonumber \vspace*{2mm}\\
& \leq 2 x_j - 2 e_{jj}^\prime x_j^\prime \leq
& \vspace*{2mm}\\
& (2d + 3) \varepsilon -
  {\displaystyle
  \left( \sum_{i=1}^n x_i^{\prime 2} - \sum_{i=1}^n x_i^2
  \right) + \sum_{i=1}^{j-1} 2 e_{ji}^\prime x_i^\prime +
  1 - \sum_{i=1}^j e_{ji}^{\prime 2}}
& \nonumber
\end{eqnarray}
for any $j \in \{ 1, 2, \ldots, n \}$.

Since
$| x_i^\prime | \leq \| f(x) \| \leq \| x \| + \varepsilon
 < d + 1$ and by Theorem \ref{thm:ndlemma}, we get
\begin{eqnarray*}
-2 (d + 1) \sum_{i=1}^{j-1} c_{ji} \varepsilon
\leq \sum_{i=1}^{j-1} 2 e_{ji}^\prime x_i^\prime
\leq 2 (d + 1) \sum_{i=1}^{j-1} c_{ji} \varepsilon.
\end{eqnarray*}
Moreover, by (\ref{eq:el1}), we have
\begin{eqnarray*}
-3 \varepsilon \leq -2 \varepsilon - \varepsilon^2
\leq 1 - \sum_{i=1}^j e_{ji}^{\prime 2}
\leq 2 \varepsilon - \varepsilon^2 \leq 3 \varepsilon.
\end{eqnarray*}
Therefore, it follows from (\ref{eq:2021-1-6an}) and
(\ref{eq:2021-1-21}) that
\begin{eqnarray*}
-\left( \left( 4 + 2 \sum_{i=1}^{j-1} c_{ji} \right) d +
        7 + 2 \sum_{i=1}^{j-1} c_{ji}
 \right) \varepsilon
\leq 2 x_j - 2 e_{jj}^\prime x_j^\prime
\leq \left( \left( 4 + 2 \sum_{i=1}^{j-1} c_{ji} \right) d +
            7 + 2 \sum_{i=1}^{j-1} c_{ji}
     \right) \varepsilon
\end{eqnarray*}
for all $j \in \{ 1, 2, \ldots, n \}$.

We note that $| x_j^\prime | < d + 1$ and
$-c_{jj} \varepsilon \leq 1 - e_{jj}^\prime
 \leq c_{jj} \varepsilon$ by Theorem \ref{thm:ndlemma}, and
since
$x_j - e_{jj}^\prime x_j^\prime =
 \big( x_j - x_j^\prime \big) +
 \big( 1 - e_{jj}^\prime \big) x_j^\prime$, we can see that
\begin{eqnarray} \label{eq:x1-x1'n}
\left| x_j - x_j^\prime \right|
\leq \left( \left( 2 + \sum_{i=1}^j c_{ji} \right) d +
            4 + \sum_{i=1}^j c_{ji}
     \right) \varepsilon
\end{eqnarray}
for $j \in \{ 1, 2, \ldots, n \}$.

Since we can select an isometry $U : D \to \mathbb{R}^n$
defined by $U(x) = x = (x_1, x_2, \ldots, x_n)$, we see that
\begin{eqnarray*}
\begin{split}
\| f(x) - U(x) \|
&    = \left\| \big( x_1^\prime - x_1, x_2^\prime - x_2,
                     \ldots, x_n^\prime - x_n
               \big)
       \right\|
  =    \left( \sum_{j=1}^n \big( x_j^\prime - x_j \big)^2
       \right)^{1/2} \\
& \leq \left( \sum_{j=1}^n
              \left( \left( 2 + \sum_{i=1}^j c_{ji} \right) d +
                     4 + \sum_{i=1}^j c_{ji}
              \right)^2
       \right)^{1/2} \varepsilon
\end{split}
\end{eqnarray*}
for all $x \in D$.
\hspace*{\fill}$\Box$
\vspace*{3mm}
%%%%%%%%%%%%%%%%%%%%%%%%%%%%%%%%%%
%%%   End of Theorem 4.1   %%%%%%%
%%%%%%%%%%%%%%%%%%%%%%%%%%%%%%%%%%

We remark that for any $i, j \in \{ 1, 2, \ldots, n \}$ with
$i \geq j$, each $c_{ij}$ is independent of $\varepsilon$ for
any `sufficiently' small $\varepsilon > 0$.

%%%%%%%%%%%%%%%%%%%%%%%%%%%%%%%%%%%%
%%%%%      5. Examples     %%%%%%%%%
%%%%%%%%%%%%%%%%%%%%%%%%%%%%%%%%%%%%

\section{Examples}

%%%%%%%%%%%%%%%%%%%%%%%%%%%%%%%%%%%%%%%
%%%%%     Example 5.1     %%%%%%%%%%%%%
%%%%%%%%%%%%%%%%%%%%%%%%%%%%%%%%%%%%%%%

\begin{exam} \label{ex:1}
We assume that $n = 4$.
We will compute some constants $c_{ij}$ by using the recurrence
formulas $(\ref{eq:diffeq1p})$ and $(\ref{eq:diffeq2p})$.
First, since $c_{11} = 1$ by $(b.1)$ in the proof of Theorem
\ref{thm:ndlemma}, it follows from $(\ref{eq:diffeq2p})$ with
$k = 1$ and $0 < \varepsilon < \frac{1}{12}$ that
\begin{eqnarray*}
c_{21}
\geq \frac{6}{11} \left( 4 + 2 \sqrt{2} + \frac{1}{12} \right)
\geq \frac{1}{2 (1 - c_{11} \varepsilon)}
     \left( 4 + 2 \sqrt{2} + \varepsilon \right),
\end{eqnarray*}
and we can choose $c_{21} = 4$ as the smallest positive integer
satisfying the last inequality.

Now, we use $(\ref{eq:diffeq1p})$ with $j = 2$ to obtain
\begin{eqnarray*}
(3c_{22} - 1) (c_{22} - 1) \geq c_{21}^2
\end{eqnarray*}
and thus, we select $c_{22} = 3$.
We note that $c_{22} = 3$ is larger than the estimate in the
proof of Theorem \ref{thm:ndlemma}.
This difference is due to the use of the formula
$(\ref{eq:diffeq1p})$ instead of $(\ref{eq:2021-2-28})$.

By $(\ref{eq:diffeq2p})$ with $k = 2$, we get
\begin{eqnarray*}
c_{32}
\geq \frac{2}{3} \left( 4 + 2 \sqrt{2} + \frac{33}{12} \right)
\geq \frac{1}{2 (1 - c_{22} \varepsilon)}
     \left( 4 + 2 \sqrt{2} + \varepsilon + 2c_{21}^2 \varepsilon
     \right)
\end{eqnarray*}
and hence, we choose $c_{32} = 7$.
Further, it follows from $(\ref{eq:diffeq1p})$ with $j = 3$ that
\begin{eqnarray*}
(3c_{33} - 1) (c_{33} - 1) \geq c_{21}^2 + c_{32}^2 = 65
\end{eqnarray*}
and hence, $c_{33} = 6$.
In the proof of Theorem \ref{thm:ndlemma}, we estimated
$c_{32} = 6$ and $c_{33} = 4$, but in this example, we estimate
the larger values for them because we use the formula
$(\ref{eq:diffeq1p})$ instead of $(\ref{eq:2021-2-28})$.

As $n = 4$, we see that
$\sigma \geq \sum\limits_{i=1}^3 c_{(i+1)i}^2
 > c_{21}^2 + c_{32}^2 = 65$ and hence,
$\varepsilon < \frac{1}{\sigma} < \frac{1}{65}$.
Therefore, we have
\begin{eqnarray*}
\frac{1}{2 (1 - c_{33} \varepsilon)}
\left( 4 + 2 \sqrt{2} + \varepsilon +
       2 \sum_{i=1}^{2} c_{3i} c_{4i} \varepsilon
\right)
< \frac{65}{118}
  \left( 4 + 2 \sqrt{2} + \frac{1}{65} + 2 \right)
< 4.88
\end{eqnarray*}
and thus, using $(\ref{eq:diffeq2p})$ with $k = 3$, we can
select $c_{43} = 5$.
In view of $(\ref{eq:diffeq1p})$ with $j = 4$, we get
\begin{eqnarray*}
(3c_{44} - 1) (c_{44} - 1)
\geq c_{21}^2 + c_{32}^2 + c_{43}^2 = 90
\end{eqnarray*}
and hence, $c_{44} = 7$, which comply with the claims of Remark
\ref{rem:3.1}.
\end{exam}

%%%%%%%%%%%%%%%%%%%%%%%%%%%%%%%%%%%%%%%
%%%%%     Example 5.2     %%%%%%%%%%%%%
%%%%%%%%%%%%%%%%%%%%%%%%%%%%%%%%%%%%%%%

\begin{exam} \label{ex:2}
We assume that $n = 5$.
As in Example \ref{ex:1} we get the following constants:
$c_{11} = 1$, $c_{21} = 4$, $c_{22} = 3$, $c_{31} = 4$,
$c_{32} = 7$, $c_{33} = 6$, $c_{41} = 4$, $c_{42} = 7$,
$c_{43} = 5$, and $c_{44} = 7$.

As $n = 5$, we see that
$\sigma \geq \sum\limits_{i=1}^4 c_{(i+1)i}^2
 > c_{21}^2 + c_{32}^2 + c_{43}^2 = 90$ and hence,
$\varepsilon < \frac{1}{\sigma} < \frac{1}{90}$.
Thus, we have
\begin{eqnarray*}
\frac{1}{2 (1 - c_{44} \varepsilon)}
\left( 4 + 2 \sqrt{2} + \varepsilon +
       2 \sum_{i=1}^3 c_{4i} c_{5i} \varepsilon
\right)
< \frac{45}{83}
  \left( 4 + 2 \sqrt{2} + \frac{1}{90} + 2 \right)
< 4.80
\end{eqnarray*}
and hence, using $(\ref{eq:diffeq2p})$ with $k = 4$, we can
select $c_{54} = 5$.
In view of $(\ref{eq:diffeq1p})$ with $j = 5$, we get
\begin{eqnarray*}
(3c_{55} - 1) (c_{55} - 1)
\geq c_{21}^2 + c_{32}^2 + c_{43}^2 + c_{54}^2  = 115
~~~\mbox{and~hence}~~~ c_{55} = 7.
\end{eqnarray*}
Moreover, due to $(\ref{eq:2021-1-23})$, we have $c_{51} = 4$,
$c_{52} = 7$, and $c_{53} = 5$.
\end{exam}

%%%%%%%%%%%%%%%%%%%%%%%%%%%%%%%%%%%%
%%%%%     Example 5.3     %%%%%%%%%%
%%%%%%%%%%%%%%%%%%%%%%%%%%%%%%%%%%%%

\begin{exam} \label{ex:3}
Let $D = \{ x \in \mathbb{R}^4 : \| x \| \leq d \}$ for some
$d \geq 1$ and let $f : D \to \mathbb{R}^4$ be a function
satisfying $f(0) = 0$ and the inequality $(\ref{eq:2021-1-13})$
for all $x, y \in D$ and some constant $\varepsilon$ with
$0 < \varepsilon < \frac{1}{90}$.
Using Theorem \ref{thm:ndtheorem} and Example \ref{ex:1}, we
can prove that there exists an isometry $U : D \to \mathbb{R}^4$
satisfying
\begin{eqnarray*}
\begin{split}
\| f(x) - U(x) \|
& \leq \left[ \sum_{i=1}^4
              \left( \left( 2 + \sum_{j=1}^i c_{ij} \right) d +
                     4 + \sum_{j=1}^i c_{ij}
              \right)^2 \,
       \right]^{1/2} \varepsilon \\
& =    \sqrt{1076 d^2 + 2376 d + 1316}\, \varepsilon \\
& <    (33 d + 37) \varepsilon
\end{split}
\end{eqnarray*}
for all $x \in D$.
\end{exam}

%%%%%%%%%%%%%%%%%%%%%%%%%%%%%%%%%%%%
%%%%%     Example 5.4     %%%%%%%%%%
%%%%%%%%%%%%%%%%%%%%%%%%%%%%%%%%%%%%

\begin{exam} \label{ex:4}
Let $D = \{ x \in \mathbb{R}^5 : \| x \| \leq d \}$ for some
$d \geq 1$ and let $f : D \to \mathbb{R}^5$ be a function
satisfying $f(0) = 0$ and the inequality $(\ref{eq:2021-1-13})$
for all $x, y \in D$ and some constant $\varepsilon$ with
$0 < \varepsilon < \frac{1}{115}$.
Using Theorem \ref{thm:ndtheorem} and Example \ref{ex:2}, we
can prove that there exists an isometry $U : D \to \mathbb{R}^5$
satisfying
\begin{eqnarray*}
\begin{split}
\| f(x) - U(x) \|
& \leq \left[ \sum_{i=1}^5
              \left( \left( 2 + \sum_{j=1}^i c_{ij} \right) d +
                     4 + \sum_{j=1}^i c_{ij}
              \right)^2 \,
       \right]^{1/2} \varepsilon \\
& =    \sqrt{1976 d^2 + 4296 d + 2340}\, \varepsilon \\
& <    (45 d + 49) \varepsilon
\end{split}
\end{eqnarray*}
for all $x \in D$.
\end{exam}

%%%%%%%%%%%%%%%%%%%%%%%%%%%%%%%%%%%%%%
%%%%%     6. Discussion     %%%%%%%%%%
%%%%%%%%%%%%%%%%%%%%%%%%%%%%%%%%%%%%%%

\section{Discussion}

We expect the Hyers-Ulam stability of isometries defined on
bounded domains to have widespread application, but no
remarkable results have been published, except for
V\"{a}is\"{a}l\"{a}'s paper \cite{va}, for 40 years after
Fickett's theorem was published:
A finite sequence $( u_0, u_1, \ldots, u_m )$ of a compact
subset $D$ of $\mathbb{R}^n$ is called a maximal sequence in
$D$ provided that $h_k := d(u_k, A_{k-1})$ is maximal in $D$
for each $k \in \{ 1, 2, \ldots, m \}$, where $A_{k-1}$ is the
affine subspace (or flat) including
$\{ u_0, u_1, \ldots, u_{k-1} \}$.
Now, among the results proved by V\"{a}is\"{a}l\"{a}, the part
related to the subject of this paper is introduced in the
following theorem.

%%%%%%%%%%%%%%%%%%%%%%%%%%%%%%%%%%%%%%
%%%%%     Theorem 6.1     %%%%%%%%%%%%
%%%%%%%%%%%%%%%%%%%%%%%%%%%%%%%%%%%%%%

\begin{thm}[V\"{a}is\"{a}l\"{a}] \label{thm:6.1}
Let $D$ be a compact subset of $\mathbb{R}^n$.
If there exist a maximal sequence $( u_0, u_1, \ldots, u_n )$
in $D$ and a constant $K \geq 1$ such that
\begin{itemize}
\item[$(i)$]  $| u_k - u_0 | \leq K h_k$ for every
              $k \in \{ 2, 3, \ldots, n \}$;
\item[$(ii)$] $D \setminus \{ u_1, \ldots, u_n \}
               \subset B_{K h_n}(u_0)$,
\end{itemize}
then for every $\varepsilon$-isometry $f : D \to \mathbb{R}^n$,
there exist an isometry $U : \mathbb{R}^n \to \mathbb{R}^n$
and a constant $K^\ast > 0$ such that
\begin{eqnarray*}
\| f(x) - U(x) \| \leq K^\ast \varepsilon
\end{eqnarray*}
for all $x \in D$.
\end{thm}

Now we have made great progress in improving Fickett's theorem
in practice using intuitive method in this paper.
As we can see, our main theorem (Theorem \ref{thm:ndtheorem})
looks much simpler and clearer than Theorem \ref{thm:6.1}
proved by V\"{a}is\"{a}l\"{a}.
According to Theorem \ref{thm:ndtheorem} of this paper,
if a function $f : D \to \mathbb{R}^n$ satisfies $f(0) = 0$ as
well as the inequality $(\ref{eq:2021-1-13})$ for all
$x, y \in D$ and for some sufficiently small constant
$\varepsilon > 0$, then there exist an isometry
$U : D \to \mathbb{R}^n$ and a constant $K > 0$ such that the
inequality $\| f(x) - U(x) \| \leq K \varepsilon$ holds for
all $x \in D$.
However, it is impossible to deduce this useful conclusion by
using Fickett's theorem.
From this point of view, we can say that Fickett's theorem has
been remarkably improved in this paper.
\vspace{5mm}

\iffalse

%%%%%%%%%%%%%%%%%%%%%%%%%%
%%%   Acknowledgment   %%%
%%%%%%%%%%%%%%%%%%%%%%%%%%

\noindent
{\bf Acknowledgments.}
This work was supported by the National Research Foundation of
Korea (NRF) grant funded by the Korea government (MSIT)
(No. 2020R1F1A1A01049560).
This work was supported by 2021 Hongik University Research Fund.
\vspace{5mm}

\fi

\end{document}